\ifx\documentclass\undefined
\documentstyle[12pt]{article}
\else
\documentclass[12pt]{article}
\usepackage{latexsym}
\usepackage{amsfonts}
\usepackage{amsmath}
\usepackage{amssymb}
\usepackage{graphicx}
\fi

\author{ D\'aniel Bezdek  \and   K\'aroly Bezdek 
\thanks{Partially supported by a Natural Sciences and 
Engineering Research Council of Canada Discovery Grant.} }

\def\kkk{\null\hfill $\Box$\smallskip}
\def\E{{\mathbb E}}

\def\S{{\mathbb S}}

\def\eps{\varepsilon}

\def\DD{{\bf D}}

\def\PP{{\bf P}}
\def\BB{{\bf B}}
\def\CC{{\bf C}}

\def\ReRe{{\bf Re}}

\newtheorem{theorem}{Theorem}[section]
\newtheorem{lemma}[theorem]{Lemma}
\newtheorem{sublemma}[theorem]{Sublemma}

\newtheorem{con}[theorem]{Conjecture}

\newcommand{\proof}{{\noindent\bf Proof:{\ \ }}}

\newcommand{\dist}{\mathop{\rm dist\,}}

\title{Shortest billiard trajectories
}

\begin{document}

\maketitle

\begin{abstract}
In this paper we prove that any convex body of the $d$-dimensional Euclidean space ($d\ge 2$) possesses at least one shortest generalized billiard trajectory moreover, any of its shortest generalized billiard trajectories is of period at most $d+1$. Actually, in the Euclidean plane we improve this theorem as follows. A disk-polygon with parameter $r>0$ is simply the intersection of finitely many (closed) circular disks of radii $r$, called generating disks, having some interior point in common in the Euclidean plane. Also, we say that a disk-polygon with parameter $r>0$ is a fat disk-polygon if the pairwise distances between the centers of its generating disks are at most $r$. We prove that any of the shortest generalized billiard trajectories of an arbitrary fat disk-polygon is a $2$-periodic one. Also, we give a proof of the analogue result for $\eps$-rounded disk-polygons obtained from fat disk-polygons by rounding them off using circular disks of radii $\eps>0$. Our theorems give partial answers to the very recent question raised by S. Zelditch on characterizing convex bodies whose shortest periodic billiard trajectories are of period $2$. 

\end{abstract}

\section{Introduction}

Mathematical billiards have been studied intensively for quite some time and have a fascinating history. (See for example the book \cite{T} of Tabachnikov for an elegant introduction to billiards.) The main motivation for writing this paper was to get an understanding of the geometry of the shortest billiard trajectories in convex bodies. The more exact details of the results proved in this paper are as follows.

\begin{figure}
	\centering
		\includegraphics [width=0.5\textwidth]{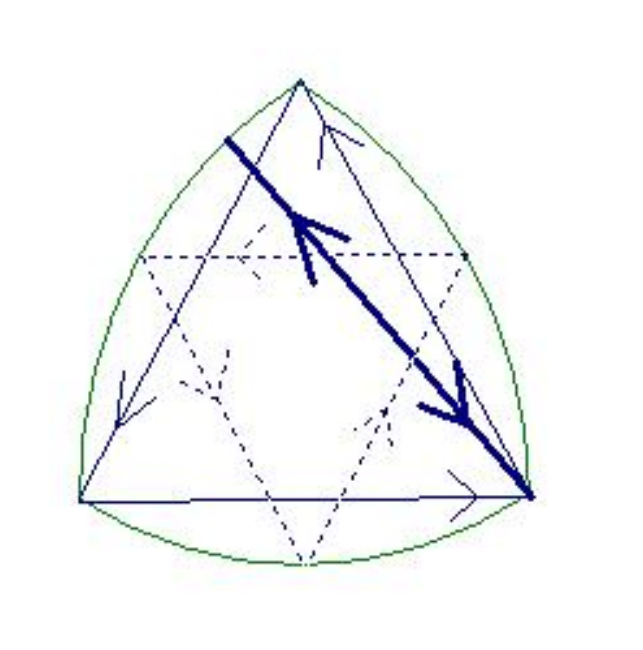}
		\caption[]{The Reuleaux triangle and some of its generalized billiard trajectories.}
	\label{fig:Reuleaux}
\end{figure}

For our purposes it seems natural to define billiard trajectories in the following way. This introduces a larger class of polygonal pathes for billiard trajectories than the traditional definition widely used in the literature (see \cite{T}). Let $\CC$ be an arbitrary convex body in the $d$-dimensional Euclidean space $\E^d, d\ge 2$ that is a compact convex set with non-empty interior in $\E^d$. Then we say that the closed polygonal path $\PP$ (possibly with self-intersections) is a {\it generalized billiard trajectory} of $\CC$ if all the vertices of $\PP$ lie on the boundary of $\CC$ and if each of the inner angle bisectors of $\PP$ between two consecutive sides of $\PP$ is perpendicular to a supporting hyperplane of $\CC$ passing through the corresponding vertex of $\PP$ (Figure~\ref{fig:Reuleaux}). If $\PP$ has $n$ sides, then we say that $\PP$ is an {\it $n$-periodic generalized billiard trajectory} in $\CC$. Note that our definition of generalized billiard trajectories coincides with the traditional definition of billiard trajectories whenever the billiard table has no corner points. It seems that the paper \cite{BC} was among the first suggesting a detailed study of generalized billiard trajectories in convex domains. For some analogue higher dimensional investigations we refer the interested reader to the recent paper \cite{G}. Generalized billiard trajectories have the following fundamental property.

\begin{theorem}\label{zero}
Let $\CC$ be a convex body in $\E^d, d\ge 2$. Then $\CC$ possesses at least one shortest generalized billiard trajectory moreover, any of the shortest generalized billiard trajectories in $\CC$ is of period at most $d+1$. 
\end{theorem}

For the sake of completeness we mention that according to the main result of \cite{BG} any $d$-dimensional billiard table with a smooth boundary, but not necessarily convex, has a $k$-periodic billiard trajectory with $k\le d+1$, which is a closely related result. Actually, in $\E^2$ we can improve Theorem~\ref{zero} as follows. In order to state that theorem in a possibly short form let us introduce the following terminology. A {\it disk-polygon} is simply the intersection of finitely many congruent closed circular disks having some interior point in common in $\E^2$.
Disk-polygons form a large class of convex domains (that is of compact convex sets with non-empty interior) in $\E^2$ including for example the well-known family of {\it Reuleaux polygons}. (See \cite{BF} for a substantial list of references and results on Reuleaux polygons and related topics.) Disk-polygons and their higher dimensional analogues which are called ball-polyhedra have many interesting geometric properties and so, it is somewhat surprising that the systematic study of this class of convex sets which is in a way a "curved analogue" of the class of convex polygons (resp., convex polytopes), has been proposed only very recently in \cite{BLNP}. In this paper we assume that whenever we take a disk-polygon, the disks generating it (simply called {\it generating disks}) are all needed that is each of them contributes to the boundary of the disk-polygon through a circular arc. Finally, let $\DD$ be a disk-polygon in $\E^2$ having the property that the pairwise distances between the centers of its generating disks of radii $r$ are at most $r$. In short, we say that $\DD$ is a {\it fat disk-polygon} with parameter $r>0$. (We note that for the purpose of this introduction it would be sufficient to assume that $r=1$. However, as in the following proofs (fat) disk-polygons are showing up with different generating radii $r$, therefore in order to accomodate them properly, we had to insist on not to fix the value of the parameter $r$.)

\begin{theorem}\label{first}
Let $\DD$ be a fat disk-polygon in $\E^2$. Then any of the shortest generalized billiard trajectories in $\DD$ is a $2$-periodic one.  
\end{theorem}

In the recent paper \cite{G}, the following fundamental question is studied that was raised by Zelditch in \cite{Z} motivated by applications to inverse spectral problems: In which convex bodies are the shortest periodic billiard trajectories of period $2$? It is proved in \cite{G} that any convex body whose inscribed ball touches the boundary of the given convex body at two diametrically opposite points has that property. Theorem~\ref{first} shows that the family of fat disk-polygons possesses that property as well.

According to Birkhoff's well-known theorem (\cite{T}) if $\BB$ is a strictly convex billiard table with smooth boundary (that is if the boundary of $\BB$ is a simple, closed, smooth and strictly convex curve) in $\E^2$, then for every positive integer $N>1$ there exist (at least two) $N$-periodic billiard trajectories in $\BB$. (In fact, here the rotation number of the billiard trajectory in question can be pre-assigned as well. Also, it is well-known that neither the convexity nor smoothness can be removed from the assumptions in order to have the same conclusion. Last but not least, for a higher dimensional analogue of the Birkhoff's theorem we refer the interested reader to \cite{FT}.) Billiard tables suitable for Birkhoff's theorem can be easily constructed from disk-polygons as follows.
Take a disk-polygon $\DD$ in $\E^2$. Then choose a positive $\eps$ not larger than the inradius of $\DD$ (which is the radius of the largest circular disk contained in $\DD$) and take the union of all circular disks of radius $\eps$ that lie in $\DD$. The set obtained in this way we call the {\it $\eps$-rounded disk-polygon} of $\DD$ and denote it by $\DD(\eps)$. The following theorem shows also the complexity of the problem of Zelditch \cite{Z} on characterizing convex domains whose shortest periodic billiard trajectories are of period $2$.

\begin{theorem}\label{second}

Let $\DD$ be a fat disk-polygon in $\E^2$. Then any of the shortest (generalized) billiard trajectories in the $\eps$-rounded disk-polygon $\DD(\eps)$ is a $2$-periodic one for all $\eps>0$ being sufficiently small.    
\end{theorem}

Actually, we believe that the following even stronger statement holds.

\begin{con}\label{con}
Let $\DD$ be a fat disk-polygon in $\E^2$. Then any of the shortest (generalized) billiard trajectories in the $\eps$-rounded disk-polygon $\DD(\eps)$ is a $2$-periodic one for all $\eps$ being at most as large as the inradius of $\DD$. 
\end{con}

\medskip

\section{Proof of Theorem~\ref{zero} }

We start with the following rather natural statement that can be proved easily with the help of Helly's theorem \cite{BF}. 

\begin{lemma}\label{lemma1}
Let $\bf F$ be a finite set of at least d+1 points and $\CC$ be a convex set in $\E^d, d\ge 2$. Then $\CC$ has a translate that covers $\bf F$ if and only if every d+1 points of $\bf F$ can be covered by a translate of $\CC$. 
\end{lemma}

\proof
For each point $p\in {\bf F}$ let $\CC_p$ denote the set of all translation vectors in $\E^d$ with which one can translate $\CC$ such that it contains $p$ i.e. let $\CC_p:=\{ t \in \E^d\ |\ p\in t+\CC\}$. Now, it is easy to see that $\CC_p$ is a convex set of $\E^d$ for all $p\in {\bf F}$ moreover, ${\bf F}\subset t+\CC$ if and only if $t\in\cap_{p\in{\bf F}}\CC_p$. Thus, Helly's theorem \cite{BF} applied to the convex sets $\{ \CC_p\ |\  p\in {\bf F}\}$ implies that ${\bf F}\subset t+\CC$ if and only if
$\CC_{p_1}\cap\CC_{p_2}\cap\dots\cap\CC_{p_{d+1}}\ne\emptyset$ holds for any $p_1, p_2,\dots ,p_{d+1}\in{\bf F}$, i.e. if and only if any $p_1, p_2,\dots ,p_{d+1}\in{\bf F}$ can be covered by a translate of $\CC$, finishing the proof of Lemma~\ref{lemma1}.
\kkk

Also the following statement is going to play a central role in our investigations. This is a generalization of the analogue $2$-dimensional statement proved in \cite{BC}.

\begin{lemma}\label{lemma2}
Let ${\bf F}=\{f_1, f_2, \dots , f_n\}$ be a finite set of points and $\CC$ be a convex body in $\E^d, d\ge 2$. Then $\bf F$ cannot be translated into the interior of $\CC$ if and only if the following two conditions hold:
There are closed supporting halfspaces $H_{i_1}^+, H_{i_2}^+, \dots , H_{i_s}^+$ of $\CC$ assigned to some points of ${\bf F}$ say, to $f_{i_1}, f_{i_2}, \dots ,  f_{i_s}$ with $1\le i_1< i_2< \dots < i_s\le n $ and a translation vector $t\in \E^d$ such that

\item (i) the translated point $t+f_{i_j}$ belongs to the closed halfspace $H_{i_j}^-$ for all $1\le j\le s$, where the interior of $H_{i_j}^-$ is disjoint from the interior of $H_{i_j}^+$ and its boundary hyperplane is identical to the boundary hyperplane of $H_{i_j}^+$ (which is in fact, a supporting hyperplane of $\CC$);

\item (ii) the intersection $\cap_{j=1}^s H_{i_j}^+ $ is nearly bounded meaning that it lies between two parallel hyperplanes of $\E^d$.
\end{lemma}

\proof
First, we assume that there are closed supporting halfspaces $H_{i_1}^+,$ $H_{i_2}^+, \dots ,H_{i_s}^+$ of $\CC$ assigned to some points of ${\bf F}$ say, to $f_{i_1}, f_{i_2}, \dots ,  f_{i_s}$ with $1\le i_1< i_2< \dots < i_s\le n $ and a translation vector $t\in \E^d$ satisfying $(i)$ as well as $(ii)$. Based on this our goal is to show that  $\bf F$ cannot be translated into the interior of $\CC$ or equivalently that $\bf F$ cannot be covered by a translate of the interior ${\rm int}\CC$ of $\CC$.
We prove this in an indirect way i.e. we assume that $\bf F$ can be covered by a translate of ${\rm int}\CC$ and look for a contradiction. Indeed, if $\bf F$ can be covered by a translate of ${\rm int}\CC$, then $t+{\bf F}$ can be covered by a translate of ${\rm int}\CC$ i.e. there is a translation vector $t^*\in \E^d$ such that
$t+{\bf F}\subset t^*+{\rm int}\CC$. In particular, if ${\bf F}^*:=\{f_{i_1}, f_{i_2}, \dots ,  f_{i_s} \}$, then $t+{\bf F}^*\subset t^*+{\rm int}\CC$. Clearly, this implies that $\cap_{j=1}^s H_{i_j}^+\subset {\rm int}\big(\cap_{j=1}^s t^*+H_{i_j}^+\big)$, a contradiction to $(ii)$. 

Second, we assume that $\bf F$ cannot be translated into the interior of $\CC$ and look for closed supporting halfspaces $H_{i_1}^+, H_{i_2}^+, \dots , H_{i_s}^+$ of $\CC$ assigned to some points of ${\bf F}$ say, to $f_{i_1}, f_{i_2}, \dots ,  f_{i_s}$ with $1\le i_1< i_2< \dots < i_s\le n $ and a translation vector $t\in \E^d$ satisfying $(i)$ as well as $(ii)$. In order to simplify matters let us start to investigate the case when $\CC$ is a smooth convex body in $\E^d$ i.e. through each boundary point of $\CC$ there exists precisely one supporting hyperplane of $\CC$. (Also, without loss of generality we can assume that the origin $o$ of $\E^d$ is an interior point of $\CC$.)    As $\bf F$ cannot be translated into ${\rm int}\CC$ therefore Lemma~\ref{lemma1} implies that there are $m\le d+1$ points of $\bf F$ say, ${\bf F}_{m}:=\{f_{j_1}, f_{j_2}, \dots ,  f_{j_m}\}$ with $1\le j_1< j_2< \dots < j_m\le n $ such that ${\bf F}_{m}$ cannot be translated into ${\rm int}\CC$. Now, let $\lambda_0:=\inf \{\lambda >0\ |\ \lambda {\bf F}_{m}\ {\rm cannot\  be\  translated \ into}\ {\rm int}\CC\}$. Clearly, $\lambda_0\le 1$ and $\lambda_0 {\bf F}_{m}$ cannot be translated into ${\rm int}\CC$ moreover, as $\lambda_0=\sup \{\delta >0 \ |\ \delta {\bf F}_{m}\ {\rm can \ be \ translated \ into\ } \CC\}$, therefore there exists a translation vector $t\in \E^d$ such that $t+\lambda_0 {\bf F}_{m}\subset \CC$. Let $t+\lambda_0 f_{i_1}, t+\lambda_0 f_{i_2}, \dots , t+\lambda_0 f_{i_s}$ with $1\le i_1< i_2< \dots < i_s\le n , 2\le s\le m \le d+1$ denote the points of $t+\lambda_0 {\bf F}_{m}$ that are boundary points of $\CC$ and let $H_{i_1}^+, H_{i_2}^+, \dots , H_{i_s}^+$ be the corresponding closed supporting halfspaces of $\CC$. We claim that ${\bf H}^+:=\cap_{k=1}^s H_{i_k}^+ $ is nearly bounded. Indeed, if ${\bf H}^+$ were not nearly bounded, then there would be a translation vector $t'\in\E^d$ with ${\bf H}^+\subset t'+{\rm int}{\bf H}^+$. As $\CC$ is a smooth convex body therefore this would imply the existence of a sufficiently small $\mu>0$ with the property that $\{ t+\lambda_0 f_{i_1}, t+\lambda_0 f_{i_2}, \dots ,  t+\lambda_0 f_{i_s}\}\subset \mu t'+{\rm int}\CC$, a contradiction. Thus, as $o\in {\rm int}\CC$ therefore the points $ f_{i_1}, f_{i_2}, \dots , f_{i_s}$ and the closed supporting halfspaces $H_{i_1}^+, H_{i_2}^+, \dots , H_{i_s}^+$ and the translation vector $t\in \E^d$ satisfy $(i)$ as well as $(ii)$. We are left with the case when $\CC$ is not necessarily a smooth convex body in $\E^d$. In this case let $\CC_N, N=1,2,\dots$ be a sequence of smooth convex bodies lying in ${\rm int}\CC$ with $\lim_{N\to +\infty}\CC_N=\CC$. As $\bf F$ cannot be translated into the interior of $\CC_N$ for all $N=1,2,\dots$ therefore applying the method described above to each $\CC_N$ and  taking proper subsequences if necessary we end up with some points of $\bf F$ say, $ f_{i_1}, f_{i_2}, \dots , f_{i_s}$ with $1\le i_1< i_2< \dots < i_s\le n$ and with $s$ convergent sequences of closed supporting halfspaces $H_{N, i_1}^+, H_{N, i_2}^+, \dots , H_{N, i_s}^+$ of
$\CC_N$ and a convergent sequence of translation vectors $t_N$ that satisfy $(i)$ and $(ii)$ for each $N$. By taking the limits $H_{i_1}^+:=\lim_{N\to +\infty}H_{N, i_1}^+,
H_{i_2}^+:=\lim_{N\to +\infty}H_{N, i_2}^+, \dots , H_{i_s}^+:=\lim_{N\to +\infty}H_{N, i_s}^+$ and $t:=\lim_{N\to +\infty}t_N$ we get the desired nearly bounded family of closed supporting halfspaces of $\CC$ and the translation vector $t\in \E^d$ satisfying $(i)$ as well as $(ii)$. This completes the proof of Lemma~\ref{lemma2}. 
\kkk

\begin{lemma}\label{lemma3}
Let $\CC$ be a convex body in $\E^d, d\ge 2$. If $\PP$ is a generalized billiard trajectory in $\CC$, then $\PP$ cannot be translated into the interior of $\CC$. 
\end{lemma}

\proof
Let $p_1, p_2, \dots , p_n$ be the vertices of $\PP$ and let $v_1, v_2, \dots , v_n$ be the points of the unit sphere $\S^{d-1}$ centered at the origin $o$ in $\E^d$ whose position vectors are parallel to the inner angle bisectors (halflines) of $\PP$ at the vertices $p_1, p_2, \dots , p_n$ of $\PP$. Moreover, let $H_1^+, H_2^+, \dots , H_n^+$ denote the closed supporting halfspaces of $\CC$ whose boundary hyperplanes are perpendicular to the inner angle bisectors of $\PP$ at the vertices $p_1, p_2, \dots , p_n$. Based on Lemma~\ref{lemma2} in order to prove that $\PP$ cannot be translated into the interior of $\CC$ it is sufficient to show that $\cap_{i=1}^{n}H_i^+$ is nearly bounded or equivalently that $o\in {\rm conv}(\{v_1, v_2, \dots , v_n\})$, where ${\rm conv}(.)$ denotes the convex hull of the corresponding set in $\E^d$. It is easy to see that $o\in {\rm conv}(\{v_1, v_2, \dots , v_n\})$ if and only if for any hyperplane $H$ of $\E^d$ passing through $o$ and for any of the two closed halfspaces bounded by $H$ say, for $H^+$, we have that $H^+\cap {\rm conv}(\{v_1, v_2, \dots , v_n\})\ne\emptyset$. Indeed, for a given $H^+$ let $t\in\E^d$ be chosen so that $t+H^+$ is a supporting halfspace of ${\rm conv}(\{p_1, p_2, \dots , p_n\})$. Clearly, at least one vertex of $\PP$ say, $p_{i_0}$ must belong to the boundary of $t+H^+$ and therefore $v_{i_0}\in H^+\cap {\rm conv}(\{v_1, v_2, \dots , v_n\})$, finishing the proof of Lemma~\ref{lemma3}. \kkk

For the purpose of the following statement it seems natural to introduce {\it generalized $(d+1)$-gons} in $\E^d$ as closed polygonal path (possibly with self-intersections) having at most $d+1$ sides.

\begin{lemma}\label{lemma4}
Let $\CC$ be a convex body in $\E^d, d\ge 2$ and let ${\cal F}_{d+1}(\CC)$ denote the family of all generalized $(d+1)$-gons of $\E^d$ that cannot be translated into the interior of $\CC$. Then ${\cal F}_{d+1}(\CC)$ possesses a minimal length member moreover, the shortest perimeter members of ${\cal F}_{d+1}(\CC)$ are identical (up to translations) with the shortest generalized billiard trajectories of $\CC$.
\end{lemma}

\proof
If $\PP$ is an arbitrary generalized billiard trajectory of the convex body $\CC$ in $\E^d$ with vertices $p_1, p_2, \dots , p_n$, then according to Lemma~\ref{lemma3} $\PP$ cannot be translated into the interior of $\CC$. Thus, by Lemma~\ref{lemma1} $\PP$ possesses at most $d+1$ vertices say, $p_{i_1}, p_{i_2}, \dots , p _{i_{d+1}}$ with $1\le i_1\le i_2\le \dots \le i_{d+1}\le n $ such that $p_{i_1}, p_{i_2}, \dots ,  p_{i_{d+1}}$ cannot be translated into the interior of $\DD$. This implies that by connecting the consecutive points of $p_{i_1}, p_{i_2}, \dots ,  p_{i_{d+1}}$ by line segments according to their cyclic ordering the generalized $(d+1)$-gon $\PP_{d+1}$ obtained,
has length $l(\PP_{d+1})$ at most as large as the length $l(\PP)$ of $\PP$ moreover, $\PP_{d+1}$ cannot be covered by a translate of ${\rm int}\CC$ i.e. $\PP_{d+1} \in {\cal F}_{d+1}(\CC)$. Now, by looking at only those members of ${\cal F}_{d+1}(\CC)$ that lie in a $d$-dimensional ball of sufficiently large radius in $\E^d$ we get via a standard compactness argument and Lemma~\ref{lemma2} that ${\cal F}_{d+1}(\CC)$ possesses a member of minimal length say, $\Delta_{d+1}(\CC)$. As the inequalities $l(\Delta_{d+1}(\CC) )\le l(\PP_{d+1})\le l(\PP)$ hold for any generalized billiard trajectory $\PP$ of $\CC$ therefore in order to finish our proof it is sufficient to show that $\Delta_{d+1}(\CC)$ is a generalized billiard trajectory of $\CC$. Indeed, as $\Delta_{d+1}(\CC)\in {\cal F}_{d+1}(\CC)$ therefore  $\Delta_{d+1}(\CC)$ cannot be translated into ${\rm int}\CC$. Thus, the minimality of $\Delta_{d+1}(\CC)$ and Lemma~\ref{lemma2} imply that if $q_1, q_2, \dots ,q_m$ denote the vertices of $\Delta_{d+1}(\CC)$ with $m\le d+1$, then there are closed supporting halfspaces $H_1^+, H_2^+, \dots , H_m^+$ of $\CC$ whose boundary hyperplanes $H_1, H_2, \dots , H_m$ pass through the points $q_1, q_2, \dots ,q_m$ (each being a boundary point of $\CC$) and have the property that $\cap_{i=1}^{m}H_i^+$ is nearly bounded in $\E^d$. If the inner angle bisector at a vertex of $\Delta_{d+1}(\CC)$ say, at $q_i$ were not perpendicular to $H_i$, then it is easy to see via Lemma~\ref{lemma2} that one could slightly move $q_i$ along $H_i$ to a new position $q_i'$ (typically an exterior point of $\CC$ on $H_i$) such that the new generalized $(d+1)$-gon $\Delta_{d+1}'(\CC)\in {\cal F}_{d+1}(\CC)$ would have a shorter length, a contradiction. This completes the proof of Lemma~\ref{lemma4}. \kkk

Finally, notice that Theorem~\ref{zero} follows from Lemma~\ref{lemma4} in a straightforward way.

\medskip

\section{Proof of Theorem~\ref{first} }

Let $\DD$ be a fat disk-polygon in $\E^2$ and assume that the conclusion of Theorem~\ref{first} does not hold for $\DD$. Based on Theorem~\ref{zero} this means that there is a shortest generalized billiard trajectory $\PP$ in $\DD$ of period $3$. Thus, $\PP$ is a generalized billiard triangle of $\DD$ whose perimeter is the smallest among the perimeters of all generalized billiard trajectories of $\DD$. Our goal is to show that this assumption leads to a contradiction, namely we will prove that ${\rm per}(\PP)>2\cdot{\rm width}(\DD)$, where  ${\rm width}(.,.)$ stands for the {\it minimal width} of the corresponding convex set denoting the smallest distance between two parrallel supporting lines of the given convex set.(This is indeed a contradiction because the double covered line segment represented by a chord of $\DD$ of length equal to the minimal width of $\DD$ is in fact, a generalized billiard trajectory of period $2$ in $\DD$.) 
   
\begin{sublemma}\label{sublemma2}
$\PP$ is a billiard triangle in $\DD$ meaning that it is a generalized billiard triangle of $\DD$ with the additional property that each vertex of $\PP$ is a smooth boundary point of $\DD$.
\end{sublemma}

\proof   
Lemma~\ref{lemma4} and Lemma~\ref{lemma2} imply the existence of three closed supporting halfplanes of $\DD$ whose boundary lines pass through the corresponding vertices of $\PP$ and whose intersection is a triangle say, ${\bf \Delta}_1^+$. Now, imagine that a vertex say, $p$ of $\PP$ is not a smooth boundary point of $\DD$. This means that we can slightly perturb the boundary line $L_1$ of the closed supporting halfplane of $\DD$ assigned to $p$ such that we create a new closed supporting halfplane for $\DD$ with boundary line $L_2$ and obtain a new triangle ${\bf \Delta}_2^+$ into which $\PP$ is inscribed. Clearly, we can arrange all this so that the inner angle bisector of $\PP$ at $p$ is not perpendicular to $L_2$. As a result we can slightly move $p$ along $L_2$ to a new position $p'$, such that the new triangle $\PP'$ obtained from $\PP$ by replacing $p$ with $p'$ possesses a strictly smaller perimeter than $\PP$ moreover, using Lemma~\ref{lemma2} it is easy to see that $\PP'$ cannot be translated into the interior of $\DD$, and therefore via Lemma~\ref{lemma4} it leads to a contradiction. This finishes the proof of Sublemma~\ref{sublemma2}.
\kkk

Let the vertices of $\PP$ be denoted by $p_1, p_2$ and $p_3$. Based on Sublemma~\ref{sublemma2} let $\DD_i$ be the uniquely determined generating disk of $\DD$ with center $c_i$ and radius $r$ whose boundary circle $S_i$ passes through $p_i$ , $1\le i\le 3$. Clearly $\PP$ is a billiard triangle of the disk-triangle ${\rm disk}{\bf \Delta}_r(\PP)=\cap_{i=1}^3\DD_i$. Hence, if $c$ denotes the center of the incircle of $\PP$, then the points $p_i, c$ and $c_i$ are all lying on a line for each $1\le i\le 3$. Also, as $\DD\subset {\rm disk}{\bf \Delta}_r(\PP)$ therefore
$$(1) \hskip1.0cm 2\cdot{\rm width}({\rm disk}{\bf \Delta}_r(\PP))\ge 2\cdot{\rm width}(\DD).$$
Moreover, as $c_1, c_2, c_3\in \DD$ therefore
$$(2) \hskip1.0cm c_1, c_2, c_3\in {\rm disk}{\bf \Delta}_r(\PP).$$

\begin{figure}
	\centering
		\includegraphics [width=0.5\textwidth]{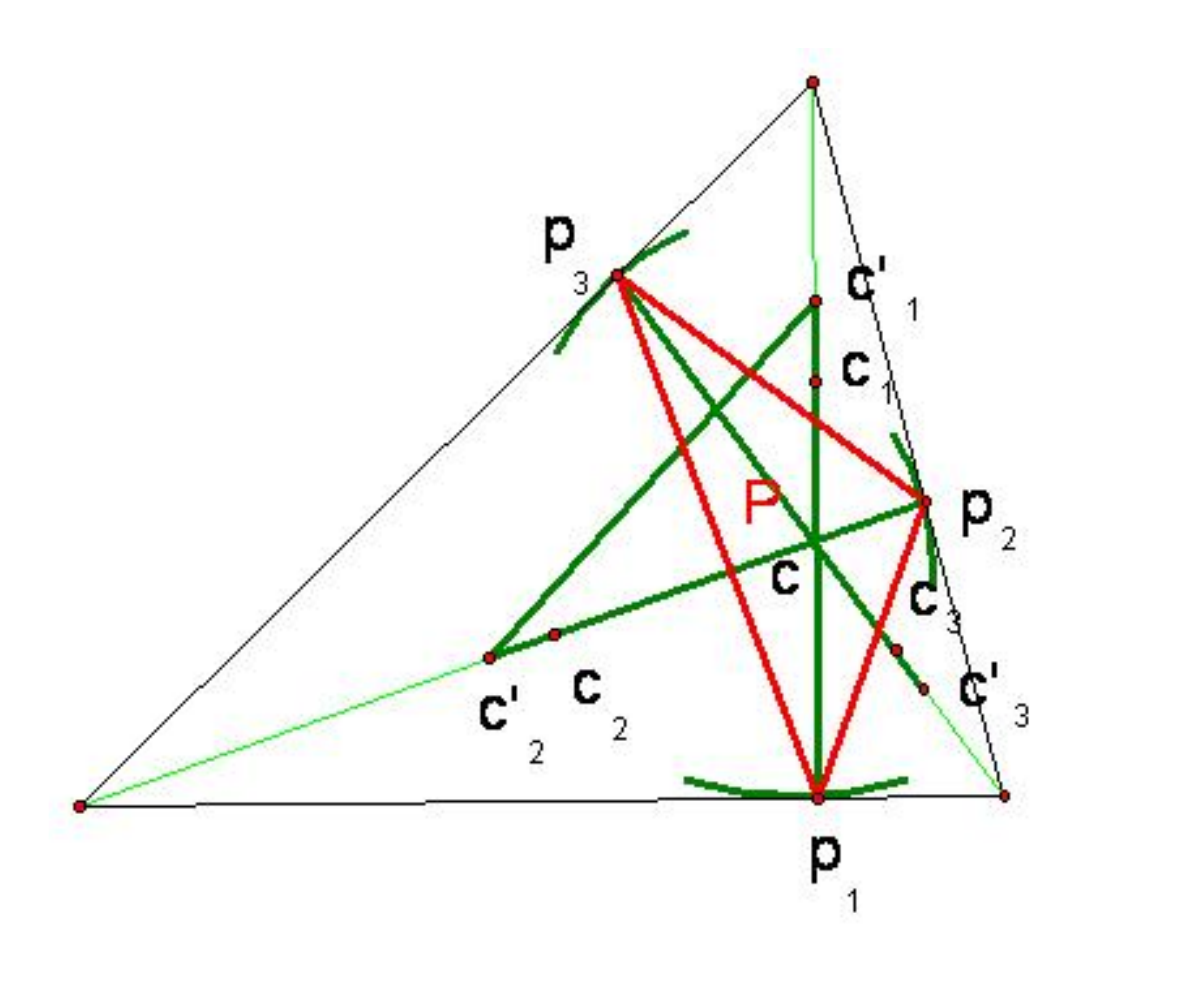}
		\caption[]{Constructing ${\rm disk}{\bf \Delta}_{r'}(\PP)$ from ${\rm disk}{\bf \Delta}_r(\PP)$.}
	\label{fig:delta'}
\end{figure}

\begin{figure}
	\centering
		\includegraphics [width=0.5\textwidth]{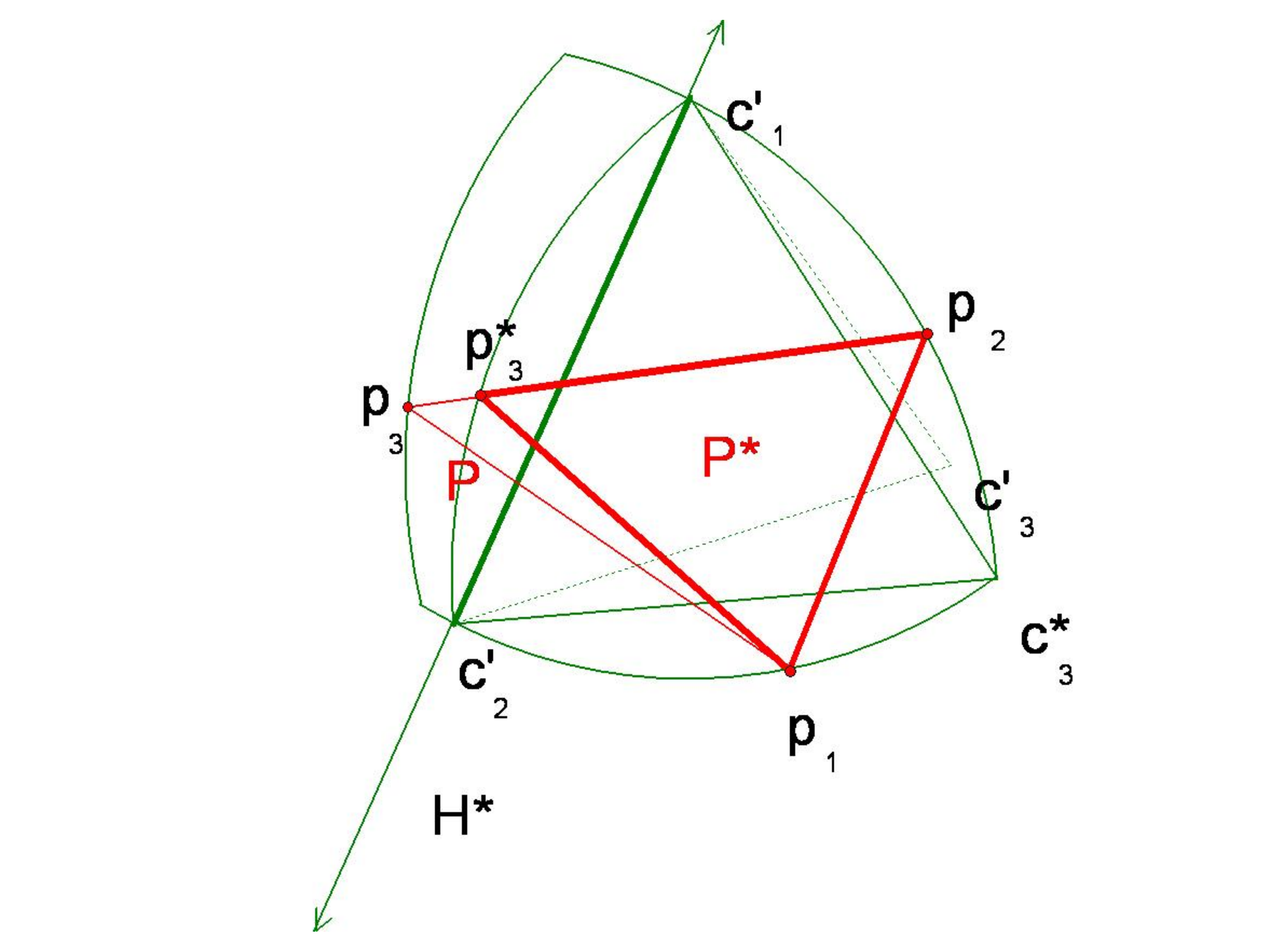}
		\caption[]{Constructing  $\PP^*$ from $\PP$. }
	\label{fig:p-star}
\end{figure}

As a next step we continuously increase the radius $r$ to a maximum value $r'$ by moving each $c_i$ along the halfline starting at $p_i$ and passing through $c$, to a new position $c'_i$ with ${\dist}(p_i, c'_i)=r'$, $1\le i\le 3$ and having the property that 
$$(2') \hskip1.0cm c'_1, c'_2, c'_3\in {\rm disk}{\bf \Delta}_{r'}(\PP),$$
where ${\dist}(. , . )$ stands for Euclidean distance between points and ${\rm disk}{\bf \Delta}_{r'}(\PP) $ denotes the intersection of the closed circular disks $\DD'_i$ of radius $r'$ centered at the points $c'_i$, $1\le i\le 3$ (Figure~\ref{fig:delta'}). It is clear that due to the maximality of $r'$ with respect to $(2')$ the diameter of the set $\{c'_1, c'_2, c'_3\}$ has to be equal to $r'$ with say, ${\dist}(c'_1, c'_2)=r'$. This means that if $S'_i$ denotes the boundary circle of $\DD'_i$, $1\le i\le 3$, then
$$(3)\hskip1.0cm c'_1\in S'_2\cap {\rm bd}({\rm disk}{\bf \Delta}_{r'}(\PP) ), c'_2\in S'_1\cap {\rm bd}({\rm disk}{\bf \Delta}_{r'}(\PP) ) , c'_3\in {\rm disk}{\bf \Delta}_{r'}(\PP),$$
where ${\rm bd}(.)$ stands for the boundary of the corresponding set. Also, as $ {\rm disk}{\bf \Delta}_r(\PP) \subset {\rm disk}{\bf \Delta}_{r'}(\PP) $ therefore 
$$(4) \hskip1.0cm 2r'=2\cdot {\rm width}( {\rm disk}{\bf \Delta}_{r'}(\PP))\ge 2\cdot {\rm width}( {\rm disk}{\bf \Delta}_{r}(\PP)).$$
Last but not least $\PP$ is a billiard triangle in the disk-triangle ${\rm disk}{\bf \Delta}_{r'}(\PP)$ as well with
$$(5) \hskip1.0cm p_i\in S'_i, 1\le i\le 3.$$
Now, let $H^*$ denote the closed halfplane bounded by the line $c'_1c'_2$ that contains $c'_3$ and let $c^*_3$ be the corner point of the $2$-sided disk-polygon $\DD'_1\cap\DD'_2$ that lies in $H^*$ (Figure~\ref{fig:p-star}).

 If $\DD^*_3$ denotes the closed circular disk of radius $r'$ centered at $c^*_3$, then $\ReRe =\DD'_1\cap\DD'_2\cap\DD^*_3$ is a Reuleaux triangle of constant width $r'$ with 
$$(6) \hskip1.0cm  \ReRe\subset {\rm disk}{\bf \Delta}_{r'}(\PP).$$
Now, $(5)$ via $(3)$ holds in the following stronger form
$$(5') \hskip1.0cm p_1\in S'_1\cap {\rm bd}(\ReRe ), p_2\in S'_2\cap {\rm bd}(\ReRe ).$$
If $S^*_3$ denotes the boundary circle of $\DD^*_3$, then $(5), (5')$ and $(6)$ imply the existence of a point $p^*_3$ with
$$(7) \hskip1.0cm p^*_3\in \PP\cap S^*_3\cap {\rm bd}(\ReRe ).$$
If $\PP^*$ denotes the triangle determined by the points $p_1, p_2$ and $p^*_3$, then Lemma~\ref{lemma2}, $(5')$ and $(7)$ imply in a straighforward way that $\PP^*$ cannot be translated into the interior of $\ReRe$.
Moreover, as $ \PP^*  \subset \PP $ therefore
$$(8)\hskip1.0cm {\rm per}(\PP )\ge {\rm per}(\PP^* ).$$

\begin{sublemma}\label{sublemma3}
${\rm per}(\PP^* )>2r'.$
\end{sublemma}

\proof
Recall the statement proved in \cite{BC} that the shortest generalized billiard trajectories of $\ReRe$ among all generalized billiard trajectories of period at most $3$ are the double covered chords of $\ReRe$ having length equal to $2r'$, where $r'$ stands for the constant width of $\ReRe$. Thus, Lemma~\ref{lemma4} implies in a straightforward way that indeed ${\rm per}(\PP^* )>2r'$. \kkk

Hence, $(1)$, $(4)$, $(8)$ and Sublemma~\ref{sublemma3} imply that ${\rm per}(\PP)>2\cdot{\rm width}(\DD)$, a contradiction. This finishes the proof of Theorem~\ref{first}.

\medskip

\section{Proof of Theorem~\ref{second} }

Assume that $\DD$ is a fat disk-polygon moreover, assume that the claim of Theorem~\ref{second} does not hold for $\DD$. In what follows, using Theorem~\ref{first}, we show that these assumptions lead to a contradiction therefore providing a proof for Theorem~\ref{second}.       

The above assumptions imply that there exists a decreasing sequence $\epsilon_n>0$ with $\lim_{n\to\infty}\epsilon_n =0$ such that for each $\epsilon_n $, the $\epsilon_n $-rounded disk-polygon $\DD(\epsilon_n )$ has a shortest (generalized) billiard trajectory $\PP(\epsilon_n )$ of period at least $3$. Now, Lemma~\ref{lemma4} implies that $\PP(\epsilon_n )$ is in fact, a (generalized) billiard triangle in $\DD(\epsilon_n )$ for all $n$. According to the well-known Blaschke selection theorem (applied to the container $\DD$ with $\DD(\epsilon_n )\subset\DD$ and to the sequence $\PP(\epsilon_n )$) one can find a subsequence of $\PP(\epsilon_n )$ whose (generalized) billiard triangles converge to a triangle or a double covered line segment labeled by $\PP$. For the sake of simplicity we assume that the subsequence in question is identical to the sequence $\PP(\epsilon_n )$ and so, we assume that the sequence $\PP(\epsilon_n )$ converges to $\PP$.  It is easy to see that $\PP$ is a generalized billiard trajectory of $\DD$ which is either a triangle or a double covered line segment. Thus, we have the following two cases:

{\it Case (A)}: $\PP$ is a generalized billiard triangle of $\DD$;
 
{\it Case (B)}: $\PP$ is a double covered line segment of $\DD$ with the property that the two lines passing through the end points of $\PP$ and being perpendicular to $\PP$ are supporting lines of $\DD$.    

{\it Case (A)}: On the one hand, clearly ${\rm per}(\PP(\epsilon_n ))\le 2\cdot {\rm width}(\DD(\epsilon_n ))$ for all $n$ and therefore 
$$(9) \hskip1.0cm {\rm per}(\PP)\le 2\cdot {\rm width}(\DD).$$ 
On the other hand, Theorem~\ref{first} and $(9)$ imply that $\PP$ must be one of the shortest generalized billiard trajectories in $\DD$ of period $2$, a contradiction. This finishes the proof of Theorem~\ref{second} in Case (A). 

{\it Case (B)}: Actually, our goal is to show that this case cannot occur. Let the vertices of $\PP(\epsilon_n )$ be denoted by $a_n, b_n$ and $c_n$ with $a_n, b_n,c_n\in{\rm bd}(\DD(\epsilon_n ))$. As $\PP(\epsilon_n )$ converges to the double covered line segment $\PP$ represented by a chord of $\DD$ therefore one can assume that for all sufficiently large $n$ two vertices of $\PP(\epsilon_n )$ say, $a_n$ and $b_n$ belong to ${\rm bd}(\DD(\epsilon_n ))$ as well as to the boundary of a closed circular disk of radius $ \epsilon_n$ centered at some interior point of $\DD$ say, $o_n$. Based on this, as $\PP(\epsilon_n )$ is a (generalized) billiard triangle in $\DD(\epsilon_n )$ therefore $\PP(\epsilon_n )$ must be an isosceles triangle  say, with ${\rm dist}(a_n, c_n)={\rm dist}(b_n, c_n)$ and with $o_n$ being the intersection of the inner angle bisectors of $\PP(\epsilon_n )$ for all sufficiently large $n$. If $2\alpha_n$ denotes the angle of the triangle $\PP(\epsilon_n )$ at the vertices $a_n$ and $b_n$, then an easy computation shows that
$$(10)\hskip0.5cm {\rm dist}(a_n, b_n)=2\epsilon_n\cos \alpha_n\ {\rm and} \ 
{\rm dist}(a_n, c_n)={\rm dist}(b_n, c_n)=\frac{\epsilon_n \cos \alpha_n}{\cos 2\alpha_n }.$$
Also, it is easy to see that
$$(11)\hskip0.5cm {\rm dist}(c_n, o_n)=\epsilon_n\sqrt{1+\frac{\cos^2 \alpha_n}{\cos^2 2\alpha_n } -2\frac{\cos^2 \alpha_n}{\cos 2\alpha_n }}.$$
Furthermore, if $d_n$ denotes the midpoint of the circular arc of radius $ \epsilon_n$ connecting $a_n$ to $b_n$ and belonging to ${\rm bd}(\DD(\epsilon_n ))$, then the lines passing through $c_n$ and $d_n$ and being perpendicular to the line $c_nd_n$ are supporting lines for $\DD(\epsilon_n )$. Thus, the double covered line segment with end points $c_n$ and $d_n$ represents a (generalized) billiard trajectory of period $2$ in $\DD(\epsilon_n )$ and so,
$$(12)\hskip0.5cm {\rm dist}(c_n, d_n)\ge \frac{1}{2} {\rm per}(\PP(\epsilon_n ) ).$$
Using $(10)$ and $(11)$ the inequality $(12)$ is equivalent to the following inequality that should hold for all sufficiently large $n$:
$$(13)\hskip0.5cm 0\ge \bigg(\cos \alpha_n + \frac{\cos \alpha_n}{\cos 2\alpha_n } \bigg)-
\bigg(1+\sqrt{1+\frac{\cos^2 \alpha_n}{\cos^2 2\alpha_n } -2\frac{\cos^2 \alpha_n}{\cos 2\alpha_n }}  \bigg).$$
Clearly, as $n$ tends to infinity, $\epsilon_n $ approaches $0$ , meanwhile the line segment with end points $c_n$ and $o_n$ converges to the chord $\PP$ of $\DD$. This forces $\alpha_n$ to converge to $\frac{\pi}{4}$.
Using Maple one can conveniently check that the function $g(x)=\bigg(\cos x + \frac{\cos x}{\cos 2x } \bigg)-
\bigg(1+\sqrt{1+\frac{\cos^2 x}{\cos^2 2x } -2\frac{\cos^2 x}{\cos 2x}}  \bigg)$ defined over the open interval $(0.7,\frac{\pi}{4})$ is a decreasing positive function. In fact, it turns out that $g(x)>0.4$ for all $0.7<x< \frac{\pi}{4}$. (Of course, this can be worked out "by hand" as well, using the proper analytic tools from calculus and thus, avoiding completely the use of Maple.) This contradicts $(13)$ finishing the proof of Theorem~\ref{second} in Case (B).

\vspace{1cm}

\medskip

\noindent
D\'aniel Bezdek,
Haskayne School of Business,
2500 University drive N.W.,
University of Calgary, Calgary, Alberta, Canada, T2N 1N4.
\newline
{\sf e-mail: daniel.bezdek@gmail.com}

\smallskip

\noindent
K\'aroly Bezdek,
Department of Mathematics and Statistics,
2500 University drive N.W.,
University of Calgary, Calgary, Alberta, Canada, T2N 1N4.
\newline
{\sf e-mail: bezdek@math.ucalgary.ca}

\end{document}